\documentclass{ifacconf}
\usepackage{cite}
\usepackage{amsmath,amssymb,amsfonts}
\newtheorem{definition}{\bf{Definition}}
\newtheorem{assumption}{Assumption}
\usepackage{graphicx}
\newtheorem{lemma}{\bf{Lemma}}
\newtheorem{remark}{Remark}
\usepackage{subcaption}
\usepackage{graphicx}
\usepackage{natbib}        
\begin{document}
\begin{frontmatter}

\title{Robust Control for a Class of Nonlinearly Coupled Hierarchical Systems with Actuator Faults} 


\author[First]{Sina Ameli} 
\author[Second]{Olugbenga Moses Anubi}

\address[First]{Department of Electrical and Computer Engineering, Florida State University (e-mail: sa19bk@my.fsu.edu).}
\address[Second]{Department of Electrical and Computer Engineering, the Center for Advanced Power Systems, Florida State University, USA (e-mail: oanubi@fsu.edu)}

\begin{abstract}                
This paper proposes an approach to addresses the control challenges posed by a fault-induced uncertainty in both the dynamics and control input effectiveness of a class of hierarchical nonlinear systems in which the high-level dynamics is nonlinearly coupled with a multi-agent low-level dynamics. The high-level dynamics has a multiplicative uncertainty in the control input effectiveness and is subjected to an exogenous disturbance input. On the other hand, the low-level system is subjected to actuator faults causing a time-varying multiplicative uncertainty in the dynamical model and associated control effectiveness. Moreover, the nonlinear coupling between the high-level and the low-level dynamics makes the problem even more challenging. To address this problem, an online parameter estimation algorithm is designed, coupled with an adaptive splitting mechanism which automatically distributes the control action among low level multi-agent systems. A nonlinear  $\mathcal{L}_2$-gain-based controller, and then a state-feedback controller are designed in the high-level, and the low-level, respectively, to recover the system from faults with high performance in the transient response, and reject the exogenous disturbance. The resulting analysis guarantees a robust tracking of the high-level reference command signal.
\end{abstract}


\end{frontmatter}

\section{Introduction}
It is well-known that the robust control can deal with any class of systems with uncertainty and disturbance (\cite{anubi2013new},~\cite{anubi2013roll}). One class of uncertainties is due to actuator faults, which causes multiplicative time-varying uncertainty in the control input matrix. In~\cite{zhang2017prescribed}, a state-feedback controller using a function such that the adaptive parameters are bounded is designed for a class of nonlinear systems with time-varying multiplicative uncertainty caused by faults in actuators. In~(\cite{stefanovski2018fault}), an $\mathcal{H}_\infty$ controller is designed in the frequency domain for LTI descriptor systems with multiplicative uncertainty due to faults and disturbances. In~(\cite{von2018stable}), an $\mathcal{H}_\infty$ controller is designed for a linear time-invariant system with disturbances as additive faults such that the $\mathcal{H}_\infty$ norm from the disturbance to the control variable is minimal, and at the same time the $\mathcal{H}_\infty$ norm from the reference to the control input is minimal. In~(\cite{hashemi2020integrated}) a robust controller using a fault estimation is designed for a class of systems with sector nonlinearity in the input subjected to exogenous signals as additive faults. The stability of the system is shown by providing sufficient conditions and the $\mathcal{L}_2$-gain performance is minimized by solving an LMI to reject the disturbance.

However, designing a robust control for nonlinear hierarchical systems, which have different levels in their structure, is challenging, and this problem is even more challenging if high-level dynamics is nonlinearly coupled with low-level dynamics. In other words the auxiliary control variable in the high-level dynamics is nonlinear itself. An effective approach to deal with hierarchical systems is designing a backstepping controller while the controller should deal with actuator faults in the low-level. In~(\cite{lan2018decoupling}), an integrated adaptive backstepping controller using a robust observer is designed for a linear time-invariant (LTI) system with disturbances including additive faults and other exogenous inputs. In~(\cite{li2019finite}), an adaptive robust backstepping controller is designed for a class of nonlinear hierarchical systems with time-varying multiplicative uncertainty due to actuator faults. The proposed controller does not need prior knowledge about the unknown terms. However, the system is linearly coupled with low-level dynamics, and the sign of the control input effectiveness needs to be known. In~(\cite{von2018stable}), a backstepping disturbance observer is designed for a fault-free nonlinear system with nonlinearly coupled hierarchical structure. However, the paper does not consider actuator faults. In~(\cite{witkowska2018adaptive}), an adaptive control allocation using backstepping approach is proposed for a nonlinear system with actuator faults, uncertainty, and disturbance. However, the backstepping control tackles with a linearly coupled hierarchical structure. Moreover, the controller can only deal with slowly-varying disturbance, and the knowledge about the fault is assumed to be known. In~(\cite{Ameli2021robust}) a $\mathcal{L}_2$-gain based controller is designed for the high-level layer and an adaptive controller is designed for the low-level layer to deal with an incipient actuator fault for a nonlinearly coupled system. In~(\cite{sassano2019optimality}), a robust optimal controller with infinite-horizon cost functional is designed for a nonlinear system with a quadratic input, which is a special case of a nonlinearly coupled hierarchical systems, and then it shows that system is $\mathcal{L}_2$-gain stable. However, this paper does not deal with faults. In~(\cite{van2018adaptive}), a robust adaptive backstepping controller is designed for a hierarchical nonlinear systems with actuator faults. However, the system is linearly coupled and the time that fault occurs is assumed to be a prior knowledge in the design.
 In~(\cite{allerhand2014robust}), the $\mathcal{L}_2$-gain analysis is provided to design a nonlinear robust controller for a linear system with uncertainty caused by actuator faults such that the controller is switched to deal with different class of uncertainties.

In this paper, a robust controller is designed for a highly nonlinear system whose high-level dynamics is subjected to disturbance, and has uncertainty, and nonlienarly coupled with low-level multi-agent systems subjected to actuator faults. In the low-level an online splitter is designed to redistribute the control law among the subsystems automatically in response to time-varying uncertainties caused by actuator faults.
Hence this paper addresses (i) the problem of nonlinear coupling between the low-level subsystems and high-level dynamics, (ii) an online redistribution of the control law for the low-level subsystems in response to actuator faults.
The remaining of the paper is organized as follows: Section 2 introduces the notation, and preliminary. Section 3, presents the problem formulation. Section 4, illustrates the control development design. Section 5, shows the numerical simulation results. Conclusion remarks are given in section 6.
\section{Notation and Preliminary}
The following notions and conventions are used throughout the paper:
$\mathbb{R}$,$\mathbb{R}^n$,$\mathbb{R}^{n\times m}$ denote the space of real numbers, real vectors of length $n$ and real matrices of $n$ rows and $m$ columns, respectively.
$\mathbb{R}_+$ denotes positive real numbers.
$X^\top$ denotes the transpose of the quantity $X$. Normal-face lower-case letters ($x\in\mathbb{R}$) are used to represent real scalars, bold-face lower-case letter ($\mathbf{x}\in\mathbb{R}^n$) represents vectors, while normal-face upper case ($X\in\mathbb{R}^{n\times m}$) represents matrices. $X>0 (\ge0)$ denotes a positive definite (semi-definite) matrix. Given appropriately dimensioned matrices $A,B,C,D$, the shorthand
\begin{align}\label{eqn:statespace}
\left[\begin{array}{c|c}A&B\\\hline C&D\end{array}\right]\triangleq C\left(sI-A\right)^{-1}B + D
\end{align}
is used to denote a state-space realization of the underlying transfer matrix.
\begin{definition}(Finite-Gain $\mathcal{L}$-stability)\cite{khalil2002nonlinear}
\noindent Consider the nonlinear system
\begin{align} \label{eqn:General}
   \mathcal{H}:\hspace{5mm} \begin{array}{rl}
    \dot{\mathbf{x}}&=f(\mathbf{x},\mathbf{v})\\
          \mathbf{z}&=h(\mathbf{x})
    \end{array}
\end{align}
where $\mathbf{x}\in\mathcal{L}_{2e}^n$, $\mathbf{v}\in\mathcal{L}_{2e}^p$, $\mathbf{z}\in\mathcal{L}_{2e}^m$ are the state, input, and output vector signals, respectively.  The system in \eqref{eqn:General}, considered as a mapping of the form $\mathcal{H}:\mathcal{L}^p_{2e}\mapsto\mathcal{L}^m_{2e}$ is said to be finite-gain $\mathcal{L}_2$-stable if there exists real non-negative constants $\gamma,\beta$ such that $ \left\|\mathcal{H}(\mathbf{v})\right\|_2\leq \gamma\left\|\mathbf{v}\right\|_2+\beta$.
\end{definition}
\begin{definition}(Dissipativity \cite{van2000l2})\label{def:dissipativity}
The dynamic system \eqref{eqn:General} is dissipative with respect to the supply rate $s(\mathbf{v},\mathbf{z})\in \mathbb{R}$, if there exists an energy function $V(\mathbf{x})\geq 0$ such that, for all $t_f\geq t_0$, 
\begin{align}\label{eqn:dissipativity ineq}
    V(\mathbf{x}(t_f))\leq V(\mathbf{x}(t_0))+\int_{t_0}^{t_f} s(\mathbf{v},\mathbf{z}) dt\hspace{2mm}\text{for all}\hspace{2mm} \mathbf{v}\in\mathcal{L}_{2e}.
\end{align}
Moreover, given a positive scalar $\gamma$, if the supply rate is taken as $s(\mathbf{v},\mathbf{z})=\gamma^2\left\|\mathbf{v}\right\|_2^2-\left\|\mathbf{z}\right\|_2^2$,
then the dissipation inequality in \eqref{eqn:dissipativity ineq} implies a finite-gain $\mathcal{L}_2$ stability and the $\mathcal{L}_2$-gain is upper bounded by $\gamma$. Consequently, the dissipativity inequality in \eqref{eqn:dissipativity ineq} becomes
\begin{align*}
    \dot{V}\leq \gamma^2\left\|\mathbf{v}\right\|_2^2-\left\|\mathbf{z}\right\|_2^2
\end{align*}
\end{definition}

\section{Problem Formulation}
The problem is to find a robust control law such that the overall system is $\mathcal{L}_2$-gain-stable. This problem considers uncertainty and disturbance in the high-level that is nonlinearly coupled with low-level dynamics with time-varying uncertainty due to the actuator faults. Consider the following class of hierarchical nonlinear uncertain systems:

\begin{align}\label{equ:nonliner sys}
   \mathcal{H}:\left\{\begin{array}{ll}
    \dot{z}=f(z,w)-g(z,w)\phi(y_1,...,y_n)
   \\[\dimexpr-\normalbaselineskip+15pt]
  y_i=\left[\begin{array}{c|c}A_i&\mathbf{b}_i\\\hline\\[\dimexpr-\normalbaselineskip+2pt] \mathbf{c}_i^\top&0\end{array}\right]u_i,\hspace{2mm}i=1,\hdots,n
\end{array}\right.,
\end{align}
where $z\in \mathbb{R}$, is the state of the high-level dynamic, $w\in \mathbb{R}$ is an unmeasurable exogenous signal, $f: \mathbb{R}\times\mathbb{R} \mapsto\mathbb{R}$, $g: \mathbb{R}\times\mathbb{R} \mapsto\mathbb{R}$, are unknown smooth nonlinear  functions, and $\phi:\mathbb{R}\times...\times\mathbb{R} \mapsto\mathbb{R}$ is a smooth nonlinear function coupling the high-level to the low-level dynamics, $y_i\in \mathbb{R}$, and $u_i\in\mathbb{R}$ are the $i$th ($i=\{1,...,n\}$) agent's low-level output and control input respectively, and $A_i\in \mathbb{R}^{n_i\times n_i}$ , $\mathbf{b}_i\in \mathbb{R}^{n_i}$ and $\mathbf{c}_i\in\mathbb{R}^{n_i}$ are the corresponding uncertain system, input and output matrices, respectively with $n_i$ the length of the associated state vector. The following assumed properties of the dynamics above are useful for subsequent developments.

\begin{assumption}\label{ass:low_level}
The low-level linear subsystems in \eqref{equ:nonliner sys} are linearly parameterizable as
\begin{align}\label{eqn:linearly_parameterized}
    \dot{\mathbf{x}}_i = A_i\mathbf{x}_i + \mathbf{b}_i u_i= Y_i(\mathbf{x}_i,u_i)\boldsymbol{\theta}_i,
\end{align}
where $\boldsymbol{\theta}_i\in \left[\begin{array}{cc}\boldsymbol{\theta}_{0_i}-\mathbf{d}_i & \boldsymbol{\theta}_{0_i}+\mathbf{d}_i\end{array}\right]$ is an uncertain parameter vector with $\boldsymbol{\theta}_{0_i}$ and $\mathbf{d}_i$ being known nominal and maximum deviation from nominal vectors respectively, and $Y_i(\mathbf{x}_i,u_i)$ is a measurable regressor matrix for which there exists a real number $\varepsilon >0$ and $T>0$ such that, for all $t\geq 0$,
\begin{align}\label{eqn:p.e.}
    \int_t^{t+T}{Y_i(\mathbf{x}_i(\tau),u_i(\tau))^\top Y_i(\mathbf{x}_i(\tau),u_i(\tau)) d\tau} \ge \varepsilon I.
\end{align}
Furthermore, the nominal values $A_{i_0},\mathbf{b}_{i_0}, \mathbf{c}_{i_0}$ of LTI model satisfy the following
\begin{enumerate}
    \item $A_{i_0}$ is Hurwitz and has no poles on the imaginary axis
    \item The nominal DC-gain is 1, i.e $\mathbf{c}_{i_0}^\top A_{i_0}^{-1}\mathbf{b}_{i_0} = 1$
\end{enumerate}
\end{assumption}
\begin{assumption}\label{as:matching condition}
For the operating points $\mathbf{x}_{0_i}$, there exists a scalar valued function $\varphi_i(\mathbf{x}_{0_i})$, such that the low-level dynamics satisfy the matching condition
\begin{align}\label{eqn:mathing_condition}
    A_i\mathbf{x}_{0_i} = \varphi_i(\mathbf{x}_{0_i})\mathbf{b}_i.
\end{align}
\end{assumption}
\begin{assumption}\label{as:Nstability}
The operating point $z_0\in\mathbb{R}$, $w_0\in\mathbb{R}$, $y_{i0}\in\mathbb{R}$, $i=1,...,n$ is a stable equilibrium point of the high level dynamics in~\eqref{equ:nonliner sys}. Thus,
\begin{align*}
    f(z_0,w_0)-g(z_0,w_0)\phi(y_{10},...,y_{n0})=0
\end{align*}
\end{assumption}

\section{Control Development}
\begin{figure}[h!]
    \centering
    \includegraphics[width=8cm, height=5cm]{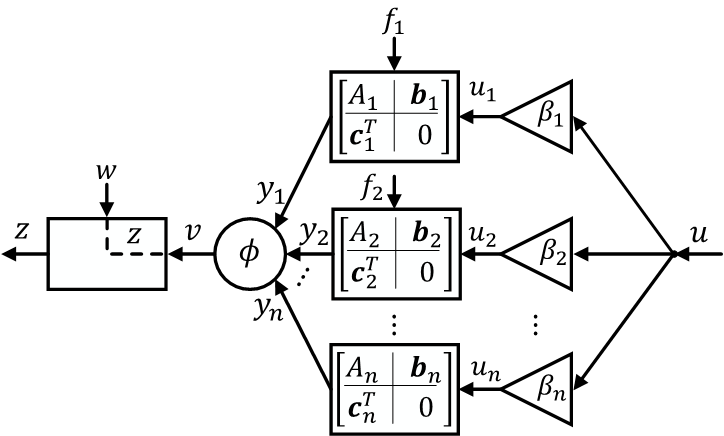}
    \caption{Hierarchical multi-agent system with nonlinear coupling.}
    \label{fig:Problem_formulation}
\end{figure}

Figure.~\ref{fig:Problem_formulation} shows a schematic of the hierarchical structure under consideration.
It illustrates a hierarchical multi-agent system with nonlinear coupling. The low-level includes multi-agent systems that each block represents a linear state-space realization subjected to a multiplicative fault changing $A_i$, and $B$. Moreover, the control input $u$ is distributed by adaptable parameter $\beta_i$ in response to the faults. Then, the agents' outputs are combined through the nonlinear function $\phi$ connecting the low-level to the high-level dynamics. The high-level dynamics is excited by the output $v$ of $\phi$. The high-level nonlinear dynamics is subject to the exogenous signal $w$. 

The high-level dynamic is actuated through the low-level dynamics coupled through the nonlinear function $\phi$. When fault occurs in any of the low-level subsystems, the corresponding parameter vector $\boldsymbol{\theta}_i$ deviates from the nominal value. The degree of the deviation corresponds to the loss of control effectiveness for that actuator. In other words, a faulty actuator can only achieve perfect tracking for smaller input. Thus, the core idea expounded upon in the control design is to track the parametric deviation using a reliable parameter estimation algorithm, then limit the corresponding control commands accordingly. Consequently, faulty agents will receive less amount of the control input while healthier ones will collectively and collaboratively compensate for this reduction in such a way that the output of the aggregator function $\phi$ remains unchanged. Thus, maintaining system response irrespective of the faulty situation. The next subsection details the parameter estimation and control allocation process.

\subsection{Parameter Estimation and Control Allocation}
The objective in this subsection is to design an algorithm that dynamically allocates the  high-level control command as references to the low-level subsystems such that healthier subsystems get more allocation and less healthy ones get less. We will refer to this algorithm as a \emph{splitter}. Since this low-level subsystems are really physics-based closed-loop actuator models, it is assumed that there are relevant parameters whose deviation from a known nominal value has a strong correlation with the health of the system. Most physical actuators have this. For instance, hydraulic actuators~(\cite{odgaard2013wind}) have natural frequencies and damping coefficients that indicate different fault conditions. Electrical actuators like batteries~(\cite{ansean2019lithium}) often have internal resistance and capacity whose values have been shown to be strong indicators of the level of degradation. Also, electric motors~(\cite{antonino2018advanced}) have internal resistance and flux parameters that are strong indication of health as well. Consequently, the internal parameters of the low-level dynamics are estimated and the resulting deviation from the respective nominal conditions are used to dynamically reallocate the high-level commands. Doing this will make the overall system automatically mitigate any faulty situation, as well as prolong the life of system by using degraded actuators less.
Consider the expanded low-level model
\begin{align}\label{eqn:aug_LTI}
    \dot{\mathbf{x}}=A\mathbf{x}+B\mathbf{u},\\
    \mathbf{y}=C\mathbf{x}
\end{align}
where $A=\textsf{blkdiag}\left(A_1,...,A_n\right)$, $B=\textsf{blkdiag}\left(\mathbf{b}_1,...,\mathbf{b}_n\right)$, $C=\textsf{blkdiag}\left(\mathbf{c}_1,...,\mathbf{c}_n\right)$, $\mathbf{u}=\begin{bmatrix}u_1&...&u_n\end{bmatrix}^\top$, $\mathbf{x}=\begin{bmatrix}\mathbf{x}_1^\top&...&\mathbf{x}_n^\top\end{bmatrix}^\top$, $\mathbf{y}=\begin{bmatrix}\mathbf{y}_1^\top&...&\mathbf{y}_n^\top\end{bmatrix}^\top$.
Convolving \eqref{eqn:aug_LTI} with the low-pass filter $H(s)=\frac{a_f}{s+a_f}$, yields
\begin{align*}
    a_f\mathbf{x}=\left(a_f I+A\right)\mathbf{x}_f+B\mathbf{u}_f,
\end{align*}
where $a_f$ is the cutoff frequency, and the filtered signals $\mathbf{x}_f$ and $\mathbf{u}_f$ are given as $\mathbf{x}_f=H(s)\mathbf{x}$ and $\mathbf{u}_f=H(s)\mathbf{u}$, respectively. Using the linearly parameterization assumption in \eqref{eqn:linearly_parameterized} yields the regression model
\begin{align*}
    a_f\left(\mathbf{x}-\mathbf{x}_f\right)=Y\left(\mathbf{x}_f,\mathbf{u}_f\right)\boldsymbol{\theta},
\end{align*}
where $Y\left(\mathbf{x}_f,\mathbf{u}_f\right) = \textsf{blkdiag}\left(Y_1(\mathbf{x}_{1_f},u_{1_f})\hdots Y_n(\mathbf{x}_{n_f},u_{n_f})\right)$ is the corresponding regressor matrix and, $\boldsymbol{\theta} = \left[\boldsymbol{\theta}_1^\top\hdots\boldsymbol{\theta}_n^\top\right]^\top$ is the combined vector of unknown parameters where $\boldsymbol{\theta}_i\in\mathbb{R}^{n_{\theta_i}}$ is the corresponding uncertain parameter vector for the $i$th low-level agent.

In order to obtain a reliable estimate of the time-varying parameter $\boldsymbol{\theta}$, a least square estimator with exponential bounded-gain forgetting factor~(\cite{slotine1991applied}) is applied as follows:
\begin{align}\label{eqn:param estimator}
     \boldsymbol{\dot{\hat{\theta}}}&=-PY\left(\mathbf{x}_f,\mathbf{u}_f\right)^\top Y\left(\mathbf{x}_f,\mathbf{u}_f\right)\boldsymbol{\hat{\theta}} \\\nonumber&\hspace{2cm}+ a_f PY\left(\mathbf{x}_f,\mathbf{u}_f\right)^\top\left(\mathbf{x}-\mathbf{x}_f\right)\\\nonumber
     \dot{P}&=\mu P-PY\left(\mathbf{x}_f,\mathbf{u}_f\right)^\top Y\left(\mathbf{x}_f,\mathbf{u}_f\right)P,
 \end{align}
where $P$ is the estimator gain matrix (estimation covariance matrix), and $\mu$ is the forgetting factor. Moreover, $Y$ persistently exciting , from assumption~\ref{ass:low_level}, implies that the estimation error $\boldsymbol{\tilde{\theta}} = \boldsymbol{\theta} - \boldsymbol{\hat{\theta}}$ converges to zero exponentially. Thus, the time-varying uncertain parameters can be tracked reliably by the estimator in \eqref{eqn:param estimator}. 
Let  
 \begin{align}\label{eqn:fault map}
     \Check{\theta}_{i}=\frac{1}{n_{\theta_i}}\sum_{j=1}^{n_{\theta_i}}{\frac{\left|\boldsymbol{\theta}_{0_j}-\boldsymbol{\hat{\theta}}_{j}\right|}{\mathbf{d}_{ij}}}.
 \end{align}
 thus, $ \Check{\theta}_{i} \in [0\hspace{1mm} 1]$ indicates the degree of fault in the $i$th low-level agent, with $ \Check{\theta}_{i}=1$ corresponding to complete loss of control effectiveness while $\Check{\theta}_{i}=0$ indicates perfectly healthy subsystem with $100\%$ control authority available.
 
 Consequently, the splitter is designed as
\begin{align}\label{eqn:spliter}
    \boldsymbol{\beta}=\frac{1}{n}\begin{bmatrix}1-\Check{\theta}_{1}\\\vdots\\1-\Check{\theta}_{q}\\1+\frac{1}{n-q}\Sigma_{i=1}^{q}\Check{\theta}_{i}\\\vdots\\1+\frac{1}{n-q}\Sigma_{i=1}^{q}\Check{\theta}_{i}
    \end{bmatrix},
\end{align}
where $q$ is the number of faulty low-level agents.
The splitter design in~\eqref{eqn:spliter} automatically redistributes the control input $u$ to low-level subsystems such that faulty actuators are given less command while the healthier ones collaboratively picks up the slack.

\subsection{Controller design}
The control design is carried out in two phases; First, the high-level design is done to regulate the high-level dynamic around the nominal operating condition against the exogenous disturbance. Next, given the desired trajectory from the high-level design, the low-level controller is designed to asymptotically track the high-level trajectory while allocating the control authorities for each low-level subsystems using the splitter in \eqref{eqn:spliter}.
\subsubsection{High-level design}
To facilitate the high-level design, the following high-level tracking error is defined
\begin{align*}
    \tilde{z}=z-z_0,
\end{align*}
where $z_0$ is a constant operating point. A corresponding filtered error is then given by:
\begin{align}\label{eqn:filtered_error}
    \rho=\tilde{z}+\eta\underbrace{\int_0^t\tilde{z} d\tau}_{\text{$\tilde{z}_I$}},
\end{align}
where $\eta>0$. Taking the derivative of~\eqref{eqn:filtered_error}, then adding and subtracting
\begin{align*}
    f(z_0,w_0)-g(z_0,w_0)\phi(y_{10},...,y_{n0}) = 0,
\end{align*} yields

\begin{align*}
    \dot{\rho}&=f(z,w)-g(z,w)\phi(y_1,...,y_n)+\eta\tilde{z},\\
              &=f(z,w)-f(z_0,w_0) - \left(l(z,w,\mathbf{y}) - l(z_0,w_0,\mathbf{y}_0)\right) +\eta \tilde{z},
\end{align*}
where
\begin{align*}
    l(z,w,\mathbf{y})\triangleq&g(z,w)\phi(y_1,...,y_n).
\end{align*}
where $y_{i0}$ is the operating point,using Assumption~\ref{as:Nstability} and invoking the mean value theorem~(\cite{rudin1964principles}), the following is obtained: 
\begin{align}\label{eqn:errordyn_high}
    \dot{\rho}=\left(\eta+h_z\right)\tilde{z}+l_w\widetilde{w}+\mathbf{l}_y^\top C\widetilde{{\mathbf{x}}},
\end{align}
where 
\begin{align*}
    h_z=&\frac{\partial f(\eta_z,\eta_w)}{\partial z}+\frac{\partial l(\eta_z,\eta_w,\boldsymbol{\eta}_y)}{\partial z},\\
    l_w=&\frac{\partial f(\eta_z,\eta_w)}{\partial w}+\frac{\partial l(\eta_z,\eta_w,\boldsymbol{\eta}_y)}{\partial w},\\
    \mathbf{l}_y=&\nabla_y l(\eta_z,\eta_w,\boldsymbol{\eta}_y),
\end{align*}
with $\eta_z = tz_0+(1-t)z$, $\eta_w = tw_0 + (1-t)w$, $\boldsymbol{\eta}_y = t\mathbf{y}_0 + (1-t)\mathbf{y}$ for some $t\in [0,\hspace{2mm}1]$,
and
\begin{align}\label{eqn:error}
    \widetilde{\mathbf{x}}&=\mathbf{x}-\mathbf{x}_0,\\
    \widetilde{w}&=w-w_0,
\end{align}
where $w_0$ is known nominal value of the exogenous disturbance. Usually, this value is used in the component design and rating. So, it is reasonable to assume that it is known. 

Here, the desire is to obtain an auxiliary control law for $\widetilde{\mathbf{x}}$ in~\eqref{eqn:errordyn_high} such that the error signal $\tilde{z}$ is robustly regulated for all $\widetilde{w}\in\mathcal{L}_{2e}$. This is then used as a reference for the low-level dynamics where the final control is designed to achieve asymptotic tracking performance on a faster time scale. The following properties of the open-loop dynamics in~\eqref{eqn:errordyn_high} are used in the subsequent design.
\begin{assumption} \label{as: Bounded control effectiveness}
The high-level dynamics is sufficiently smooth. Thus, the uncertain terms $h_z$, $l_w$, and $\mathbf{l}_y$ in~\eqref{eqn:errordyn_high} are bounded as follows; there exists $\mathbf{l}_0\in\mathbb{R}^n$ and positive constants $\alpha, \bar{h}_z$ such that
\begin{align}\label{eqn:conic}
    \mathbf{l}_y^\top CC^\top\mathbf{l}_0\geq\alpha\\
    \left|h_z\right|\leq \bar{h}_z,\\
    \left|l_w\right|\leq \bar{l}_w.
\end{align}

\end{assumption}

Consequently, the auxiliary control law is designed as:
\begin{align}\label{eqn:aux_ctrl}
    \widetilde{\mathbf{x}}=-k_1C^\top\mathbf{l}_0\rho,
\end{align}
where $k_1>0$ is a control gain and $\mathbf{l}_0$ satisfies the conic constraint in~\eqref{eqn:conic}. Thus, the corresponding high-level closed loop error system is given by 
\begin{align}\nonumber
    \dot{\rho}&=\left(\eta+h_z\right)\tilde{z}-k_1\mathbf{l}_{y}^\top CC^\top\mathbf{l}_0\rho+l_w\widetilde{w}\\\nonumber
              &=\left(\eta+h_z\right)\left(\rho-\eta \tilde{z}_I\right)-k_1\mathbf{l}_{y}^\top CC^\top\mathbf{l}_0\rho+l_w\widetilde{w}\\\label{eqn:cls-high}
              &=\left(\eta+h_z-k_1\mathbf{l}_{y}^\top CC^\top\mathbf{l}_0\right)\rho-\eta\left(\eta+h_z\right)\tilde{z}_I+l_w\widetilde{w}.
\end{align}
The following theorem gives the robust performance of the high-level auxiliary control law in~\eqref{eqn:aux_ctrl}.

\begin{thm}
Consider the high-level auxiliary control law in ~\eqref{eqn:aux_ctrl}. Given $\gamma>0$, if the control gain is chosen to satisfy the sufficient condition 
\begin{align}\label{eqn:sc1}
    k_1\ge\frac{\left(\bar{h}_z+2\eta\right)^2}{4\alpha\eta} + \frac{\bar{l}_w^2}{4\alpha\gamma^2} + \frac{1}{\alpha},
\end{align}
then the corresponding closed-loop error system in~\eqref{eqn:cls-high} is $\mathcal{L}_2$-gain stable and the $\mathcal{L}_2$-gain from the exogenous disturbance $\widetilde{w}$ to the regulation error $\tilde{z}$ is upper bounded by $\gamma$.
\end{thm}

\begin{pf}
Consider the energy function
\begin{align*}
V=\frac{1}{2}\rho^2+\frac{1}{2}\eta^2\tilde{z}_I^2.
\end{align*}
Since $\left\|\rho\right\|_2 = \left\|\left(1+\frac{\eta}{s}\right)\tilde{z}\right\|_2\ge\left\|\tilde{z}\right\|_2$, in accordance with the definition and results in Definition~\ref{def:dissipativity}, it suffices to show that $\dot{V}\le \gamma^2\widetilde{w}^2 - \rho^2$. Taking first time derivative of $V$ then adding and subtracting the term $\gamma^2\widetilde{w}^2-\rho^2$ yields
\begin{align*}
    \dot{V}&=\rho\dot{\rho} + \eta^2z_I\dot{z}_I\\
           &=\left(\eta+h_z-k_1\mathbf{l}_{y}^\top CC^\top\mathbf{l}_0\right)\rho^2  -\eta\left(\eta+h_z\right)\rho z_I + l_w\rho\widetilde{w}\\
           &\hspace{1cm} +\eta^2\left(\rho-\eta z_I\right)z_I\\
           &=\left(\eta+h_z-k_1\mathbf{l}_{y}^\top CC^\top\mathbf{l}_0+1\right)\rho^2  -\eta h_z\rho z_I + l_w\rho\widetilde{w}\\
           &\hspace{1cm}-\eta^3\tilde{z}_I^2-\gamma^2\widetilde{w}^2+\left(\gamma^2\widetilde{w}^2-\rho^2\right)\\
           &=\left(\eta+h_z-k_1\mathbf{l}_{y}^\top CC^\top\mathbf{l}_0+1+\frac{h_z^2}{4\eta^2}+\frac{l_w^2}{4\gamma^2}\right)\rho^2\\
           &\hspace{5mm}-\eta\left(\eta z_I + \frac{h_z}{2\eta}\rho\right)^2 -\gamma^2\left(\widetilde{w}-\frac{l_w}{2\gamma^2}\rho\right)^2 + \left(\gamma^2\widetilde{w}^2-\rho^2\right)\\
           &\leq-\left(k_1\alpha-\eta-\left|{h}_z\right|-\frac{\bar{l}_w^2}{4\gamma^2}-1-\frac{{h}_z^2}{4\eta^2}\right)\rho^2 +\left(\gamma^2\widetilde{w}^2-\rho^2\right)\\
           &\leq-\left(k_1\alpha-\frac{\left(h_z+2\eta\right)^2}{4\eta}-\frac{\bar{l}_w^2}{4\gamma^2}-1\right)\rho^2 +\left(\gamma^2\widetilde{w}^2-\rho^2\right)\\
           &\leq\gamma^2\widetilde{w}^2-\rho^2.
\end{align*}
\end{pf}
In the next subsection, the low-level control law is designed to achieve asymptotic tracking of the high-level auxiliary input in~\eqref{eqn:aux_ctrl}.
\subsubsection{Low-level design}
The objective in the low level control design is to improve the tracking performance for faulty low-level systems using the splitter design in~\eqref{eqn:spliter}. 
Consider the following low-level tracking error:
\begin{align*}
    \mathbf{e}=\widetilde{\mathbf{x}}+k_1C^\top\mathbf{l}_0\rho.
\end{align*}
Taking the first time-derivative yields the low-level open-loop error system
\begin{align*}
    \dot{\mathbf{e}}&=\dot{\mathbf{x}}+k_1C^\top\mathbf{l}_0\dot{\rho}\\
                    &=A\mathbf{x} + B\mathbf{u} + k_1C^\top\mathbf{l}_0\dot{\rho}\\
                    &=A\left(\mathbf{e}+\mathbf{x}_0-k_1C^\top\mathbf{l}_0\rho\right) + B\mathbf{u} + k_1C^\top\mathbf{l}_0\dot{\rho}\\
                    &=A\mathbf{e} + B\left(\mathbf{u}+\varphi(\mathbf{x}_0)\right) + \mathbf{w}_\rho,
\end{align*}
where $\mathbf{w}_\rho \triangleq k_1\left(C^\top\mathbf{l}_0\dot{\rho}-AC^\top\mathbf{l}_0\rho\right)$ is an unknown exogenous signal, which is bounded from the high-level design. 
The low-level control, restricted to the splitter direction, is then designed as
\begin{align}\label{eqn:ll_ctrl}
    \mathbf{u}=-\varphi(\mathbf{x_0}) - \boldsymbol{\beta}\mathbf{k}_2^\top\mathbf{e},
\end{align}
where $\boldsymbol{\beta}^\top=\begin{bmatrix} \beta_1&...&\beta_n
\end{bmatrix}$, with $\displaystyle \sum_{i=1}^n\beta_i=1$, is given in \eqref{eqn:spliter}. Thus, the low-level closed-loop error system is given by:
\begin{align}\label{eqn:ll_closed}
    \dot{\mathbf{e}} = \left(A-B\boldsymbol{\beta}\mathbf{k}_2^\top\right)\mathbf{e} + \mathbf{w}_\rho.
\end{align}

The following theorem gives the sufficient condition for the control gain  $\mathbf{k}_2$ to robustly regulate the closed loop error $\mathbf{e}$.

\begin{thm}
Consider the low-level control law in \eqref{eqn:ll_ctrl}. Given $\alpha_l>0$, if the control gain $\mathbf{k}_2$ is chosen to satisfy 
\begin{align}\label{eqn:ll_suff_cond}
    \left(2\underline{a}_r+\alpha_l\right)I -  B\boldsymbol{\beta}\mathbf{k}_2^\top - \mathbf{k}_2\boldsymbol{\beta}^\top B^\top \le 0,
\end{align}
where $\underline{a}_r = \max\textsf{Re}\left(\textsf{eig}\{A\}\right)$, then the closed loop error system in \eqref{eqn:ll_closed} is finite-gain $\mathcal{L}_2$-stable and the $\mathcal{L}_2$-gain with respect to the exogenous input $\mathbf{w}_\rho$ is upper bounded by $\frac{\lambda_1}{\lambda_2}$, where $\lambda_1$ and $\lambda_2$ satisfy
\begin{align}
    \frac{1}{\lambda_1} + \lambda_2 = \alpha_l.
\end{align}
\end{thm}

\begin{pf}
Consider the energy function
\begin{align*}
    V = \frac{1}{2}\left\|\mathbf{e}\right\|^2.
\end{align*}
Taking the first time-derivative and substituting the closed loop error system in \eqref{eqn:ll_closed} yields
\begin{align*}
    \dot{V} = \frac{1}{2}\mathbf{e}^\top\left(A^\top+A - B\boldsymbol{\beta}\mathbf{k}_2^\top - \mathbf{k}_2\boldsymbol{\beta}^\top B^\top\right)\mathbf{e} + \mathbf{e}^\top\mathbf{w}_\rho,
\end{align*}
which, after using the Young's inequality and adding and subtracting the term $\frac{\lambda_2\left\|\mathbf{e}\right\|^2}{2}$, becomes
\begin{align*}
    \dot{V} &\le \frac{1}{2}\mathbf{e}^\top\left(A^\top+A-B\boldsymbol{\beta}\mathbf{k}_2^\top - \mathbf{k}_2\boldsymbol{\beta}^\top B^\top + \alpha_l I\right)\mathbf{e} \\&\hspace{2cm}+ \frac{\lambda_2}{2}\left(\frac{\lambda_1}{\lambda_2}\left\|\mathbf{w}_\rho\right\|^2-\left\|\mathbf{e}\right\|^2\right).
\end{align*}
Thus, after using the sufficient condition in \eqref{eqn:ll_suff_cond}, it follows that
\begin{align*}
    \dot{V}\le \frac{\lambda_2}{2}\left(\frac{\lambda_1}{\lambda_2}\left\|\mathbf{w}_\rho\right\|^2-\left\|\mathbf{e}\right\|^2\right),
\end{align*}
which shows that the low-level closed loop error system is finite-gain $\mathcal{L}_2$ stable with $\mathcal{L}_2$-gain upper bounded by $\frac{\lambda_1}{\lambda_2}$.
\end{pf}

\begin{remark}
Suppose $B$ can be decomposed as $B = B_0 + \Delta B$, with $\left\|\Delta B\right\|\le\sigma_b$, where $B_0$ is known and $\sigma_b>0$ is a known upper bound on the uncertainty $\Delta B$. Consider the choice 
\begin{align*}
    \mathbf{k}_2 = \varepsilon B_0\boldsymbol{\beta},
\end{align*}
where $0<\varepsilon\le2$. Then, a sufficient condition in terms of $A$ only can be specified as
\begin{align*}
    2\underline{a}_r+\alpha_l+\sigma_b\le 0.
\end{align*}
It is straightforward to obtain the inequality about by applying the Young's inequality on the term $\Delta B\boldsymbol{\beta}\boldsymbol{\beta}^\top B_0^\top$ with the parameter $\varepsilon$.
\end{remark}
\section{Numerical Simulation}
In this section, the proposed control is validated on a 5MW variable pitch wind turbine model using Fatigue, Aerodynamics, Structures, and Turbulence (FAST) simulator developed by the US national renewable energy laboratory (NREL)~(\cite{jonkman2009definition}), and the implementation of the proposed controller is available online~\footnote{https://github.com/Sina-eng/Robust-control-and-estimator-for-Wind-Turbine-with-Actuator-Fault-Simulation-FAST-NREL-}. A lumped-parameter model of the rotor dynamics obtained in (\cite{wasynczuk1981dynamic}) gives the high-level dynamics as
 \begin{align}\label{eqn:WT}
 \begin{array}{rcl}
          f(z,w)&=&\frac{cw^3}{2Jz}\left(\frac{w}{z}-m_1\right)\textsf{e}^{(-m_2\frac{w}{z})}-\frac{P_0}{Jz},\\
    g(z,w)&=&\frac{cw^3}{6Jz}m_3\textsf{e}^{(-m_2\frac{w}{z})},\\
    \phi(\mathbf{y})&=&\left\|\mathbf{y}\right\|_2^2,
\end{array}
 \end{align}
 where $z$ is the rotor speed,  $\mathbf{y}\in\mathbb{R}^3$ is a vector of the pitch angle, $w$ is the wind speed, $m_1=5.4184$, $m_2=0.0682$, and $m_3 = 0.029$ are positive constants obtained experimentally,  $P_0 = 5296610 \textsf{W}$ is the rated mechanical power, $c=9.6E5$ is a positive constant, and $J = 43784700 \textsf{ kg-m}^2$ is the total drive-train inertia. There are three actuators with the following state-space dynamics
 \begin{align}\label{eqn:actuator model}
     \dot{\mathbf{x}}_i=\begin{bmatrix}
     0&1\\
     -\omega_{ni}^2&-2\zeta\omega_{ni}
     \end{bmatrix}\mathbf{x}_i
     +\begin{bmatrix}
     0\\
     \omega_{ni}^2
     \end{bmatrix}u_i,\hspace{5mm} i=1,2,3,
 \end{align}
where $\mathbf{x}_i=\begin{bmatrix}x_{1i}\\x_{2i}\end{bmatrix}$, $x_{1i}$, and $x_{2i}$ are pitch angle and pith rate for each actuator, respectively, $\zeta_i$ is the damping ratio, and $\omega_{ni}$ is the natural frequency. The variation of the damping coefficient and natural frequency have been shown to be accurate indicator of faults in the hydraulic system~\cite{odgaard2013wind}. For the least square parameter estimation, the following regression model quantities are used:
 \begin{align*}
     \Check{{x}}_i&=a_f x_{2i}-a_f x_{2fi},\\
     \boldsymbol{\theta}_i&=\begin{bmatrix}
     \omega_{ni}^2\\
     2\zeta_i\omega_{ni}
     \end{bmatrix},\\
     Y_i&=\begin{bmatrix}
     u_{fi}-x_{1fi}&-x_{2fi}
     \end{bmatrix}.
 \end{align*}
The deviation indicator is then given as
 \begin{align*}
     \Check{\theta}_{i}=\frac{1}{2}\left(\frac{\left|\omega_{n0}^2-\widehat{\omega_n^2}_i\right|}{d_\omega}+\frac{\left|(2\zeta\omega_n)_{0}-\widehat{2\zeta\omega_n}_i\right|}{d_{\zeta}}\right),
 \end{align*}
 where $\omega_{n0}^2=123.4321$, $(2\zeta\omega_n)_0=13.332$, $d_\omega=111.7357 $, $d_\zeta=10.254$.
 Moreover, the estimator uses a bounded gain matrix to tune the forgetting factor as follows~(\cite{slotine1991applied})
 \begin{align*}
     \mu(t)=\mu_0\left(1-\frac{\left\|P\right\|}{k_0}\right),
 \end{align*}
 where $P$ is the gain matrix, $\mu_0$, $k_0$ are the maximum forgetting rate, and bound for the induced norm of the gain matrix. This techniques prevents the gain matrix $P$ from becoming unbounded in case the excitation is not strong enough.
 Also the operating rotor speed is $z_0=1.267 \textsf{ rad/sec}$. The design parameters are chosen such that the sufficient condition in~\eqref{eqn:sc1} holds. To obtain the bounds on uncertainties, the bounds $11.4 \textsf{ m/sec}\leq w\leq 25 \textsf{ m/sec}$ on the wind speed are used.  Consequently, the bounds $\bar{h}_z=2.54$, $\bar{l}_w=7.8$, $\alpha=3$, $\mathbf{l}_0=\begin{bmatrix}
 1&1&1 \end{bmatrix}^\top$ are obtained. Then, for a value of $\gamma=0.3$, the sufficient condition in~\eqref{eqn:sc1} is satisfied with the choice $k_1 = 61$. Moreover, the inequality in \eqref{eqn:ll_suff_cond} is satisfied with the choice $\mathbf{k}_2 = \left[\begin{array}{cccccc}50&1&50&1&50&1\end{array}\right]$.

\begin{figure}[h!]
    \centering
    \includegraphics[width=8.5cm,height=3.5cm]{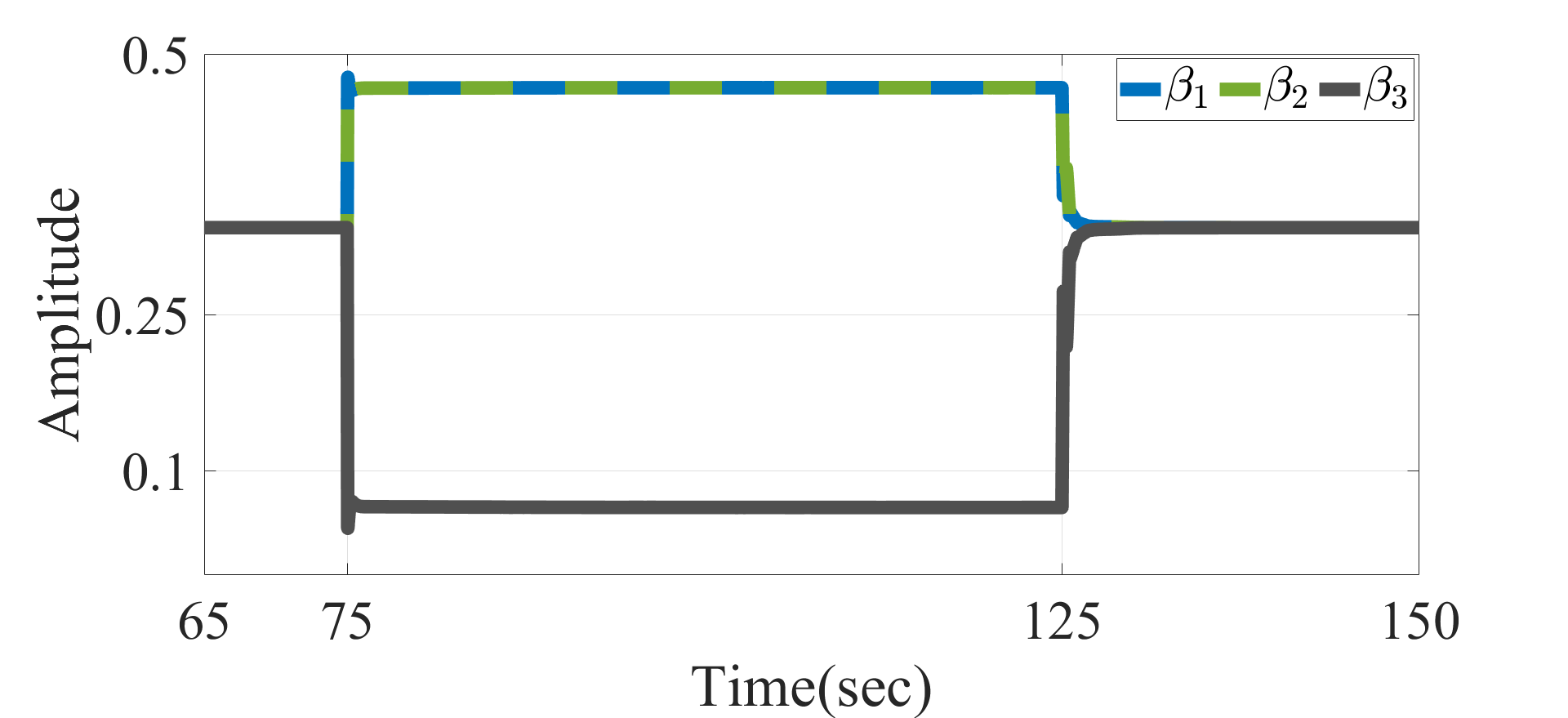}
    \caption{Automatic distribution of the control input among actuators}
    \label{fig:splitter}
\end{figure}

\begin{figure}[h!]
    \centering
    \includegraphics[width=8.5cm,height=3.5cm]{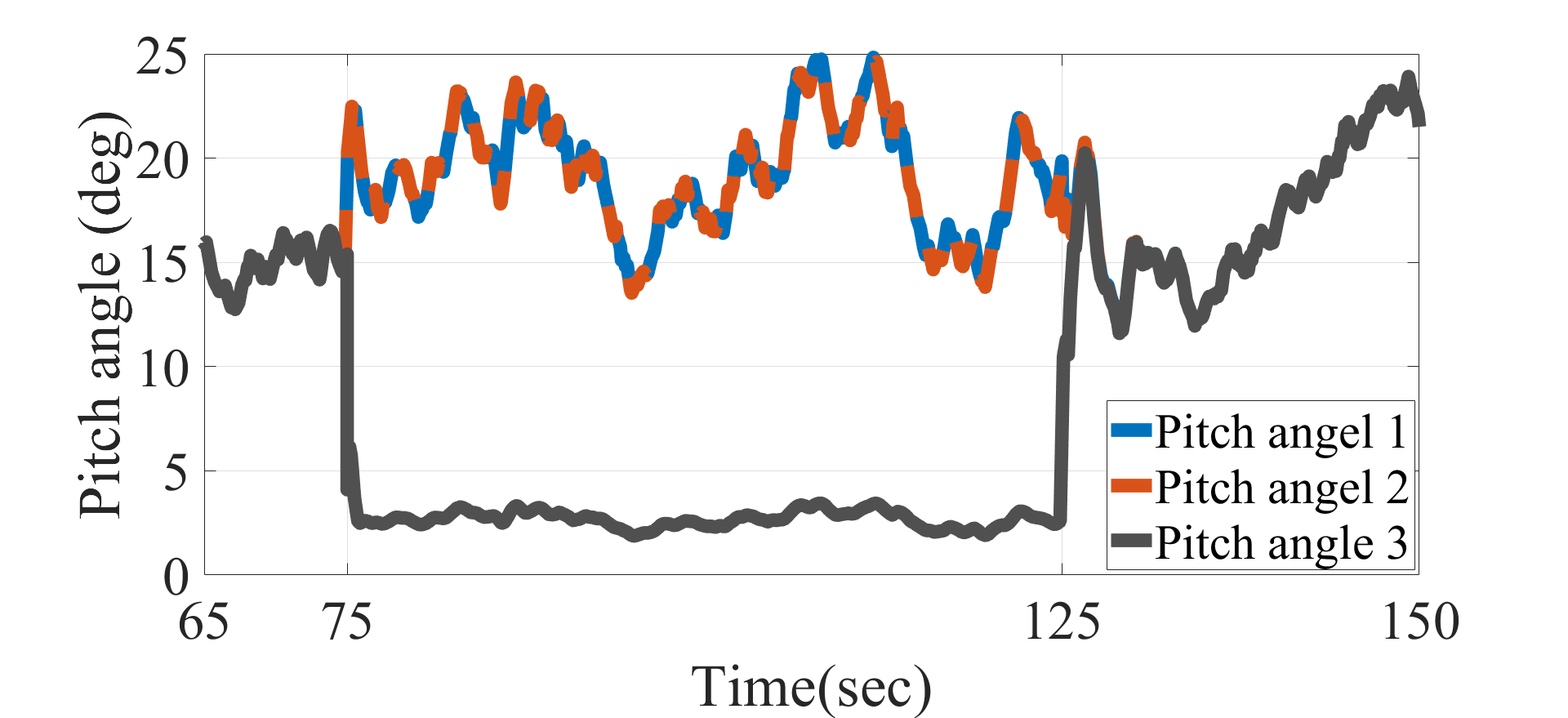}
    \caption{Pitch angle responses}
    \label{fig:pitch_sth}
\end{figure}

\begin{figure}[h!]
    \centering
    \includegraphics[width=8.5cm,height=4cm]{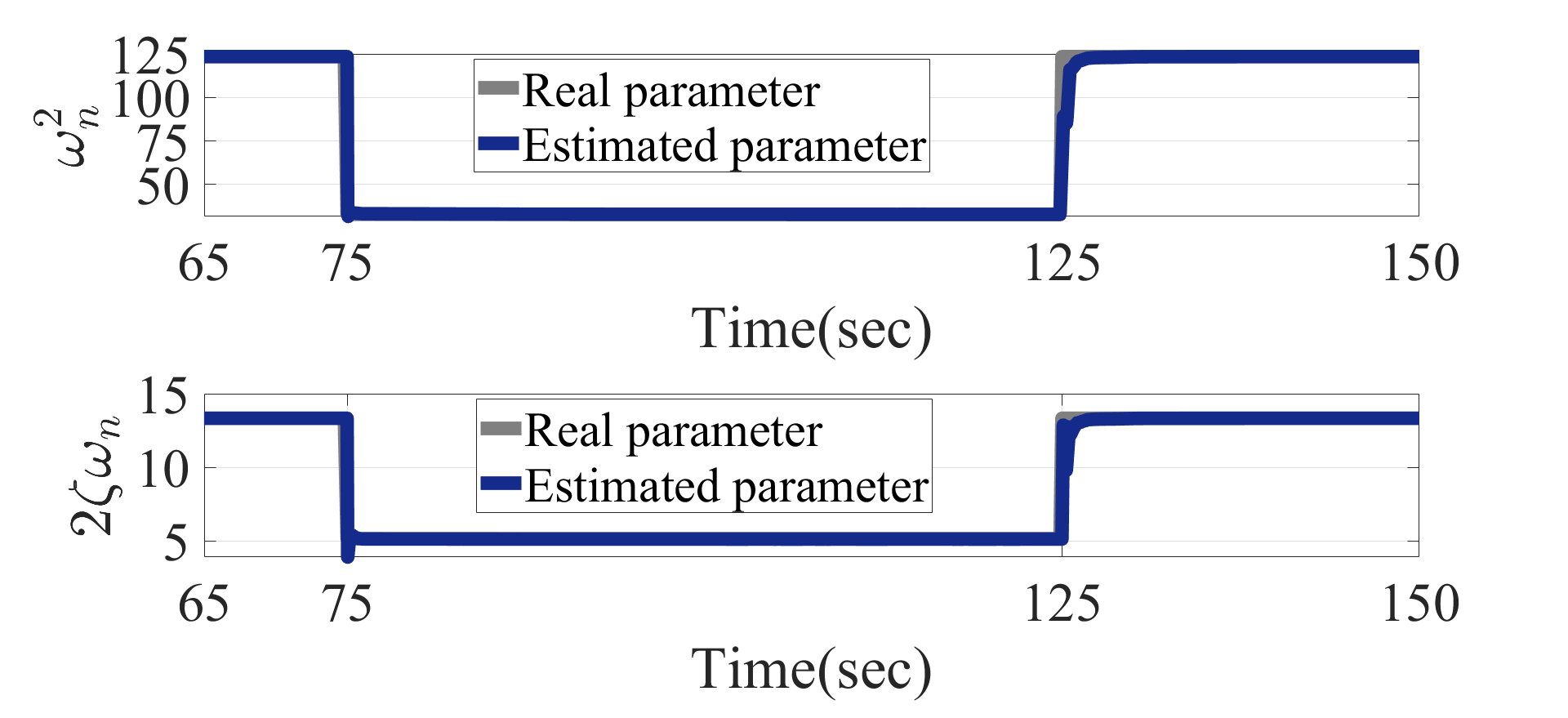}
    \caption{Online parameter identification of $\omega_n^2$, and $2\zeta\omega_n$}
    \label{fig:wn}
\end{figure}

\begin{figure}[h!]
    \centering
    \includegraphics[width=8.5cm,height=3.5cm]{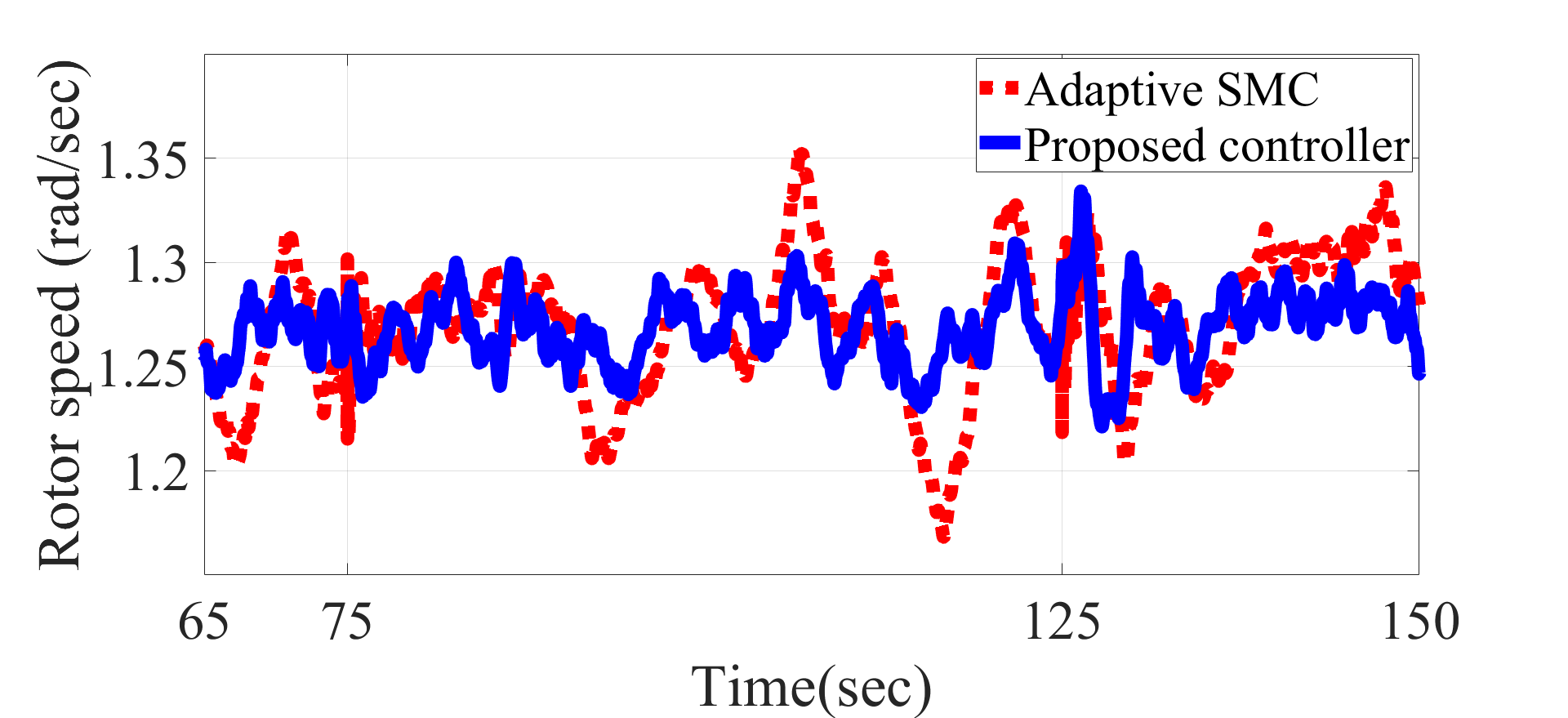}
    \caption{Rotor speed response}
    \label{fig:rotorspeed}
\end{figure}

\begin{figure}[h!]
    \centering
    \includegraphics[width=8.5cm,height=3.5cm]{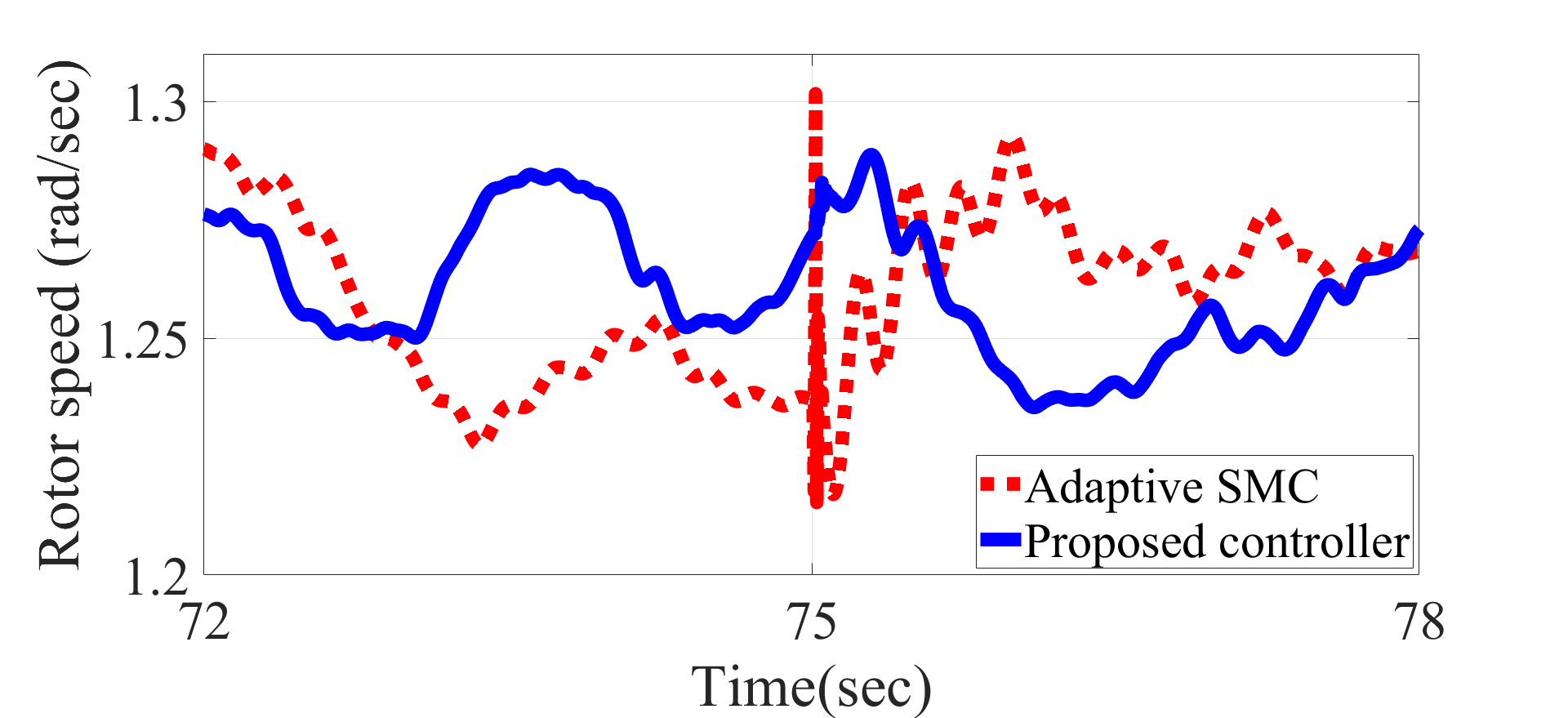}
    \caption{Rotor speed response}
    \label{fig:rotorspeed226}
\end{figure}

 The rotor speed response for the proposed controller is compared with an adaptive integral sliding mode control (SMC)~(\cite{ameli2019adaptive}). A stochastic wind profile with the mean value of $w_0 = 22 \textsf{ m/sec}$ is applied. 
 For this simulation scenario, the third actuator is faulty, with fault occurring abruptly at $75 \textsf{ sec}$, and vanishes at $125 \textsf{ sec}$. The automatic distribution of control input by the proposed splitter is shown in Fig.~\ref{fig:splitter}.
It shows that the faulty pitch actuator receives less control input in response to the fault.
Figure.~\ref{fig:pitch_sth} shows the pitch angle responses for the proposed controller. 
\begin{figure}[h!]
    \centering
    \includegraphics[width=8.5cm,height=3.5cm]{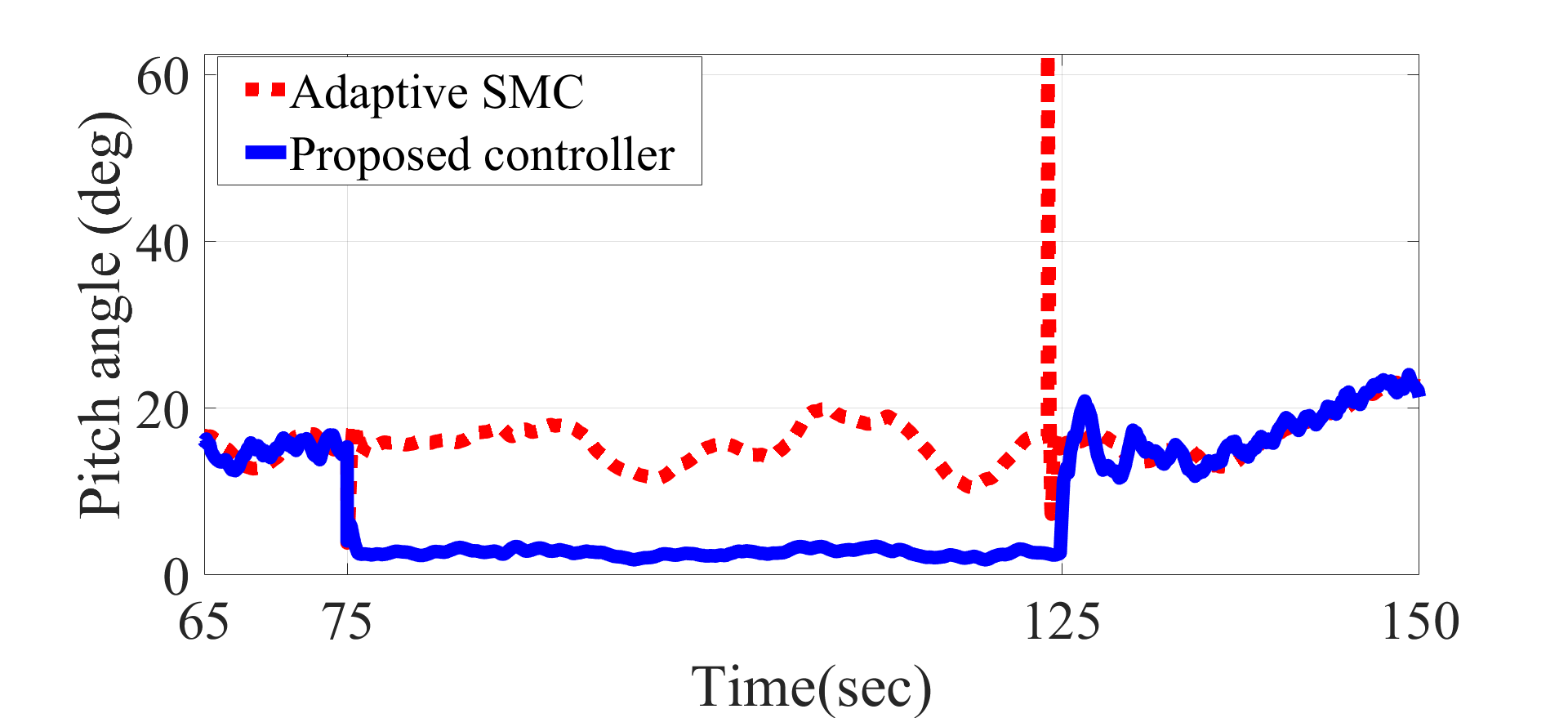}
    \caption{Faulty pitch angle response}
    \label{fig:pitch_sth_compare}
\end{figure}
It shows that when the third actuator is faulty, the other two healthy actuators are collectively collaborating to compensate for the faulty actuator. Note that both healthy actuators have the same response due to the splitter design. The faulty pitch angles of the two controller are shown in Fig.~\ref{fig:pitch_sth_compare}. It shows that the adaptive SMC has huge spike at 125 sec when the fault vanishes abruptly.  
The online parameter identification is shown in Fig.~\ref{fig:wn}.
It shows that the estimator is fast and precise in tracking the time-varying parameters.
The rotor speed response is shown in Fig.~\ref{fig:rotorspeed}, and Fig.~\ref{fig:rotorspeed226}. It shows that the proposed controller has less fluctuations especially it significantly outperforms the adaptive SMC when the fault occurs at 75 sec.

\section{Conclusion and future work}
This paper addressed the problem for a class of nonlinearly coupled hierarchical systems including multi agents subjected to actuator faults whose outputs should be collectively controlled. A splitter using parameter estimation along with a controller is proposed in response to the faults such that they collectively track a desired output required for the high-level dynamics. It was shown that the high-level closed-loop system is $\mathcal{L}_2$-gain-stable, while the error in the low-level is asymptotically stable. The results show that the splitter improves the transient response. Future works for this class of system include: 1) obtaining the minimum upper bound of the $\mathcal{L}_2$-gain of the system to improve the robustness of the controller; 2) making the robust controller less conservative 3) sensor faults, and failure.
\bibliography{ifacconf}  

                                                   







\end{document}